\documentclass[12pt]{article}
\usepackage{amsmath}
\usepackage{amssymb}
\usepackage{amsfonts}
\usepackage{amsthm}
\usepackage{blkarray}
\usepackage{arydshln}
\usepackage{multirow}
\usepackage{float}
\usepackage[dvipdfmx]{graphicx}

\newcommand{\bigzerou}{%
\smash{\lower1.7ex\hbox{\bg 0}}}
\theoremstyle{definition}
\newtheorem{theorem}{Theorem}[section]

\newtheorem{defi}{Definition}[section]
\newtheorem{cor}{Corollary}[section]
\newtheorem{conj}{Conjecture}[section]

\newtheorem{ex}{Example}[section]
\newtheorem{Rem}{Remark}[section]
\newtheorem{lem}{Lemma}[section]

\newcommand{\ba}{\begin{eqnarray}}
\newcommand{\ea}{\end{eqnarray}}
\newcommand{\ban}{\begin{eqnarray*}}
\newcommand{\ean}{\end{eqnarray*}}
\newcommand{\no}{\nonumber}

\def\d{{\partial}}

\begin{document}
\title{\begin{flushright}
\begin{minipage}[b]{5em}
\normalsize
${}$ \\
\end{minipage}
\end{flushright}
{\bf Correspondence between Gravitational Virtual Structure Constants and Gravitational Gromov-Witten Invariants of Projective Calabi-Yau Hypersurfaces}}
\author{Masao Jinzenji\\
\\
\it Department of Mathematics, Okayama University\\
\it Okayama, 700-8530, Japan\\
{\it e-mail address: pcj70e4e@okayama-u.ac.jp}}
\maketitle

\begin{abstract}
In this paper, we clarify the relation between the generalized hypergeometric functions used in the mirror computation of a Calabi-Yau hypersurface in $CP^{k-1}$ and gravitational Gromov-Witten invariants 
of a certain type. As a result, we give a proof of Givental's version of the mirror theorem for the Calabi-Yau hypersurface from the point of view of virtual structure constants. Our result also includes an extension of 
Givental's version of the mirror theorem. Moreover, we define ``ghost Gromov-Witten invariants'' $\langle
{\cal O}_{h^{k-2}}{\cal O}_{h^{-1}}\rangle_{0,d}$ and discuss their properties in the last part of this paper.  
\end{abstract}
\section{Introduction}
Let $M_{k}^{k}$ be a smooth degree $k$ hypersurface in $CP^{k-1}$, which is a Calabi-Yau manifold of 
complex dimension $k-2$. Let us introduce here the main ingredient of this paper: the gravitational 
virtual structure constants \cite{JM,JM2}, which were originally defined as intersection numbers on
$\widetilde{Mp}_{0,2}(k,d)$, the moduli space 
of quasimaps from $CP^{1}$ with two marked points to $CP^{k-1}$ of degree $d$.
\begin{defi}\footnote{In our previous papers \cite{JM, JM2}, we denoted the gravitational virtual structure constants by $w(\sigma_{j}({\cal O}_{h^{a}}){\cal O}_{h^{b}})_{0,d}$. However, that is the physicists' notation, and the 
mathematicians' notation $w(\tau_{j}({\cal O}_{h^{a}}){\cal O}_{h^{b}})_{0,d}$ differs from it by an overall factor ($\sigma_{j}=j!\cdot\tau_{j}$). Our normalization uses the mathematicians' notation, and we change the notation in this paper.
}
Let us introduce the following polynomial:
$$e(x,y)=\prod_{r=0}^{k}(rx+(k-r)y).$$
We define the gravitational virtual structure constants $w(\tau_{j}({\cal O}_{h^{a}}){\cal O}_{h^{b}})_{0,d}\;\;(d\geq 1)$ of $M_{k}^{k}$
by the following multi-variable residue integral: 
\ba
w(\tau_{j}({\cal O}_{h^{a}}){\cal O}_{h^{b}})_{0,d} &:= &\frac{1}{(2 \pi \sqrt{-1})^{d+1}} \oint_{C_{0}} \frac{dz_{0}}{(z_{0})^{k}} \oint_{C_{1}} \frac{dz_{1}}{(z_{1})^{k}} \dots \oint_{C_{d}} \frac{dz_{d}}{(z_{d})^{k}} \no\\
&& \times (z_{0})^{a} (z_{1} - z_{0})^{j} \prod_{l=1}^{d} e(z_{l-1},z_{l}) \prod_{l=1}^{d-1} \frac{1}{kz_{l} (2z_{l} - z_{l-1} - z_{l+1})} \no\\
&& \times (z_{d})^{b}.
\label{residue}
\ea
Here, $\frac{1}{(2 \pi \sqrt{-1})} \oint_{C_{0}}dz_{0}$ and $\frac{1}{(2 \pi \sqrt{-1})} \oint_{C_{d}}dz_{d}$ represent the operation of taking the residue at $z_{0}=0$ and at $z_{d}=0$, respectively, and
$\frac{1}{(2 \pi \sqrt{-1})} \oint_{C_{l}}dz_{l}\;\;(l=1,2,\cdots,d-1)$ means the operation of taking the sum of the residues at $z_{l}=0$ and $z_{l}=\frac{z_{l-1}+z_{l+1}}{2}$.
We also define the virtual structure constant $w({\cal O}_{h^{a}}{\cal O}_{h^{b}})_{0,d}$ by
\ba
w({\cal O}_{h^{a}}{\cal O}_{h^{b}})_{0,d}:=w(\tau_{0}({\cal O}_{h^{a}}){\cal O}_{h^{b}})_{0,d}.
\ea
\label{wvsdefi}
\end{defi}
\begin{Rem}
From the above definition, we can easily see that
\ban
w(\tau_{j}({\cal O}_{h^{a}}){\cal O}_{h^{b}})_{0,d}\neq 0\Longrightarrow j+a+b=k-3.
\ean 
\end{Rem}

\begin{defi}
\ba
\tilde{L}^{k,k}_{j}(x):=1+\sum_{d=1}^{\infty}\frac{d}{k}\cdot w({\cal O}_{h^{k-2-j}}{\cal O}_{h^{j-1}})_{0,d}\cdot e^{dx}.
\label{diffL}
\ea
\label{wgdef}
\end{defi}
\begin{defi}
\ba
J_{j,0}^{k,k}(x)&:=&1,\no\\
J_{j,m}^{k,k}(x)&:=&\int^{x}\tilde{L}^{k,k}_{j-m+1}(u)J_{j,m-1}^{k,k}(u)du\quad\quad(m\geq 1).
\ea
Here, the integration constant is fixed by the following rule:
\ban
\int^{x}u^mdu&:=&\frac{1}{m+1}x^{m+1},\no\\
\int^{x}u^me^{\alpha u}du&:=&\sum_{j=0}^{m}e^{\alpha x}\left(\frac{(-1)^{j}}{\alpha^{j+1}}{}_{m}P_{j}x^{m-j}\right)\quad(\alpha\neq 0)\;\;.
\label{int}
\ean
\end{defi}
\begin{theorem}
If $1\leq m\leq j$, the following equality holds:
\ba
J_{j,m}^{k,k}(x)=\frac{x^{m}}{m!}+\sum_{n=1}^{m}\frac{x^{m-n}}{(m-n)!}\left(\sum_{d=1}^{\infty}\frac{1}{k}w(\tau_{n-1}({\cal O}_{h^{k-2-j+m-n}}){\cal O}_{h^{j-m}})_{0,d}e^{dx}\right).
\ea
\label{mainB}
\end{theorem}
\begin{defi}
We define the mirror map, i.e., the coordinate change of the deformation parameter, as follows:
\ba
t=t(x)=J^{k,k}_{1,1}(x):=x+\sum_{d=1}^{\infty}\frac{1}{k}\cdot w({\cal O}_{h^{k-3}}{\cal O}_{h^{0}})_{0,d}\cdot e^{dx}.
\ea
\label{mmap}
\end{defi}
At this stage, we introduce the genus $0$ Gromov-Witten invariant $\langle\prod_{a=0}^{k-2}({\cal O}_{h^a})^{m_{a}}\rangle_{0,d}$ of degree $d$.
Here, we omit the rigorous definition of this invariant, but we mention some of its properties.
In general, it is non-zero only if the following condition is satisfied: 
\ba
&&\langle\prod_{a=0}^{k-2}({\cal O}_{h^a})^{m_{a}}\rangle_{0,d}\neq 0\;\Longrightarrow \sum_{a=0}^{k-2}m_{a}(a-1)=k-5.
\label{sel}
\ea  
If $d=0$, it vanishes except in the following case:
\ba
\langle{\cal O}_{h^a}{\cal O}_{h^{b}}{\cal O}_{h^{c}}\rangle_{0,0}=k\cdot\delta_{k-2,a+b+c}.
\ea
If $d\geq 1$, it satisfies the following equality:
\ba
\langle\prod_{a=0}^{k-2}({\cal O}_{h^a})^{m_{a}}\rangle_{0,d}=\delta_{m_0,0}d^{m_{1}}\langle\prod_{a=2}^{k-2}({\cal O}_{h^a})^{m_{a}}\rangle_{0,d}.
\ea
This follows from the {\bf Puncture Axiom}:
\ba
\langle{\cal O}_{h^{0}}\prod_{a=0}^{k-2}({\cal O}_{h^a})^{m_{a}}\rangle_{0,d}=0\;\;(d\geq 1),
\ea
and the {\bf Divisor Axiom}:
\ba
\langle{\cal O}_{h}\prod_{a=0}^{k-2}({\cal O}_{h^a})^{m_{a}}\rangle_{0,d}=d\langle\prod_{a=0}^{k-2}({\cal O}_{h^a})^{m_{a}}\rangle_{0,d}\;\;(d\geq 1),
\ea 
for Gromov-Witten invariants.

In this paper, we also consider gravitational Gromov-Witten invariants \cite{ehx} of the following type:
\ba
\langle\tau_{n}({\cal O}_{h^{p}})\prod_{a=0}^{k-2}({\cal O}_{h^a})^{m_{a}}\rangle_{0,d}\quad  (\sum_{a=0}^{k-2}m_{a}\leq 2).
\label{GOT}
\ea
It is non-zero only if the following condition is satisfied: 
\ba
&&\langle\tau_{n}({\cal O}_{h^{p}})\prod_{a=0}^{k-2}({\cal O}_{h^a})^{m_{a}}\rangle_{0,d}\neq 0\;\Longrightarrow n+p-1+\sum_{a=0}^{k-2}m_{a}(a-1)=k-5.
\label{sel}
\ea 
If $d=0$ and $n\geq 1$, it vanishes because the moduli space of genus $0$ stable curves with at most $3$ marked points is trivial. If $d\geq 1$, it obeys the {\bf Puncture Axiom}: 
\ba
\langle{\cal O}_{h^{0}}\tau_{n}({\cal O}_{h^{p}})\prod_{a=0}^{k-2}({\cal O}_{h^a})^{m_{a}}\rangle_{0,d}=\langle\tau_{n-1}({\cal O}_{h^{p}})\prod_{a=0}^{k-2}({\cal O}_{h^a})^{m_{a}}\rangle_{0,d},
\label{punc}
\ea
the {\bf Dilaton Axiom}:
\ba
\langle\tau_{1}({\cal O}_{h^{0}})\prod_{a=0}^{k-2}({\cal O}_{h^a})^{m_{a}}\rangle_{0,d}=(-2+\sum_{a=0}^{k-2}m_{a})\langle\prod_{a=0}^{k-2}({\cal O}_{h^a})^{m_{a}}\rangle_{0,d},
\label{dila}
\ea
and {\bf Hori's Equation} \cite{hori}:
\ba
\langle({\cal O}_{h})\tau_{n}({\cal O}_{h^{p}})\prod_{a=0}^{k-2}({\cal O}_{h^a})^{m_{a}}\rangle_{0,d}=d\langle\tau_{n}({\cal O}_{h^{p}})\prod_{a=0}^{k-2}({\cal O}_{h^a})^{m_{a}}\rangle_{0,d}
+\langle\tau_{n-1}({\cal O}_{h^{p+1}})\prod_{a=0}^{k-2}({\cal O}_{h^a})^{m_{a}}\rangle_{0,d}.\no\\
\label{hori}
\ea
Another well-known property is Witten's celebrated {\bf Genus $0$ Topological Recursion
Relation} \cite{witten}:
\ba
\langle\tau_{n}({\cal O}_{h^{p}}){\cal O}_{h^q}{\cal O}_{h^r}\rangle_{0,d}=
\sum_{f=0}^{d}\sum_{a=0}^{k-2}\frac{1}{k}\langle\tau_{n-1}({\cal O}_{h^{p}}){\cal O}_{h^a}\rangle_{0,f}\langle{\cal O}_{h^{k-2-a}}{\cal O}_{h^q}{\cal O}_{h^r}\rangle_{0,d-f}.
\label{witrr}
\ea
By using these properties, gravitational Gromov-Witten invariants of the type 
given in (\ref{GOT}) are determined from genus $0$ Gromov-Witten invariants.

With this setup, we state the following theorem, which is another main result of this paper. 
\begin{theorem}
By inverting the mirror map, we obtain the following equality:
\ba
J_{j,m}^{k,k}(x(t))=\frac{t^{m}}{m!}+\sum_{n=1}^{m}\frac{t^{m-n}}{(m-n)!}\left(\sum_{d=1}^{\infty}\frac{1}{k}\langle\tau_{n-1}({\cal O}_{h^{k-2-j+m-n}}){\cal O}_{h^{j-m}}\rangle_{0,d}e^{dt}\right).
\ea
\label{mainA}
\end{theorem}
\begin{Rem}
In the $m=1$ case, a geometrical proof has already been presented in \cite{Jin6}. {\bf The $m=j$ case corresponds 
to the statement of Givental's mirror theorem for the $J$-function.} We will explain the connection in Section 3.
\end{Rem}
\vspace{2em}
\noindent
{\bf Acknowledgment} 
We would like to thank Prof. K. Kuwata and Prof. K. Matsuzaka for discussions. 
Our research is partially supported by JSPS grants No. 22K03289 and No. 24H00182.

\section{Proof of Theorem \ref{mainB} and Theorem \ref{mainA}}
\begin{lem}
\ba
&&\sum_{f=1}^{d-1}(d-f)\frac{1}{k}w(\tau_{n-1}({\cal O}_{h^{k-2-j+m-n}}){\cal O}_{h^{j-m}})_{0,f}\cdot
\frac{1}{k}w({\cal O}_{h^{k-2-j+m}}{\cal O}_{h^{j-m-1}})_{0,d-f}\no\\
&=&d\cdot\frac{1}{k}w(\tau_{n}({\cal O}_{h^{k-2-j+m-n}}){\cal O}_{h^{j-m-1}})_{0,d}+
\frac{1}{k}w(\tau_{n-1}({\cal O}_{h^{k-1-j+m-n}}){\cal O}_{h^{j-m-1}})_{0,d}\no\\
&&-\frac{1}{k}w(\tau_{n-1}({\cal O}_{h^{k-2-j+m-n}}){\cal O}_{h^{j-m}})_{0,d}.
\label{wtrr}
\ea
\label{wtrl}
\end{lem}
\noindent
{\bf Proof)} Since $w(\tau_{n-1}({\cal O}_{h^{k-2-j+m-n}}){\cal O}_{h^{j-m}})_{0,f}\cdot
\frac{1}{k}w({\cal O}_{h^{k-2-j+m}}{\cal O}_{h^{j-m-1}})_{0,d-f}$ is represented by the following residue integral:
\ban
&&\frac{1}{(2 \pi \sqrt{-1})^{d+1}} \oint_{C_{0}} \frac{dz_{0}}{(z_{0})^{k}} \oint_{C_{1}} \frac{dz_{1}}{(z_{1})^{k}} \dots \oint_{C_{d}} \frac{dz_{d}}{(z_{d})^{k}} (z_{0})^{k-2-j+m-n} (z_{1} - z_{0})^{n-1} \prod_{l=1}^{d} e(z_{l-1},z_{l}) \no\\
&&\times\left(\prod_{l=1}^{f-1} \frac{1}{kz_{l} (2z_{l} - z_{l-1} - z_{l+1})}\right)\frac{1}{kz_{f}}\left(\prod_{l=f+1}^{d-1} \frac{1}{kz_{l} (2z_{l} - z_{l-1} - z_{l+1})}\right)(z_{d})^{j-m-1}\no\\
&=&\frac{1}{(2 \pi \sqrt{-1})^{d+1}} \oint_{C_{0}} \frac{dz_{0}}{(z_{0})^{k}} \oint_{C_{1}} \frac{dz_{1}}{(z_{1})^{k}} \dots \oint_{C_{d}} \frac{dz_{d}}{(z_{d})^{k}} (z_{0})^{k-2-j+m-n} (z_{1} - z_{0})^{n-1} \prod_{l=1}^{d} e(z_{l-1},z_{l}) \no\\
&&\times(2z_{f}-z_{f-1}-z_{f+1})\left(\prod_{l=1}^{d-1} \frac{1}{kz_{l} (2z_{l} - z_{l-1} - z_{l+1})}\right)(z_{d})^{j-m-1},
\ean
the assertion of the lemma follows from the equality:
\ban
\sum_{f=1}^{d-1}(d-f)(2z_{f}-z_{f-1}-z_{f+1})=d(z_{1}-z_{0})+z_{0}-z_{d}.\hspace{3cm}\Box
\ean
\begin{Rem}
By introducing the notation of multi-point gravitational virtual structure constants and Hori's equation
used in \cite{JM}, (\ref{wtrr}) is rewritten as follows:
 \ba
&&\sum_{a=0}^{k-2}\sum_{f=0}^{d}w(\tau_{n-1}({\cal O}_{h^{k-2-j+m-n}}){\cal O}_{h^{a}})_{0,f}\cdot
\frac{1}{k}w({\cal O}_{h^{k-2-a}}{\cal O}_{h^{j-m-1}}|{\cal O}_{h})_{0,d-f}\no\\
&=&w(\tau_{n}({\cal O}_{h^{k-2-j+m-n}}){\cal O}_{h^{j-m-1}}|{\cal O}_{h})_{0,d}.
\label{wtrr2}
\ea 
From this representation, we can regard (\ref{wtrr}) as the virtual version of the genus $0$ topological recursion relation for gravitational Gromov-Witten invariants. 
\end{Rem}
\vspace{1em}
\noindent
{\bf Proof of Theorem \ref{mainB})} The theorem can be straightforwardly proved by induction on $m$.
However, we write down the detailed computation here for clarity.  
In the $m=1$ case, by using the definition of $\tilde{L}^{k,k}_{j}(x)$ (\ref{diffL}), we obtain 
\ban
J_{j,1}^{k,k}(x)&=&\int^{x}du\left(1+\sum_{d=1}^{\infty}\frac{d}{k}\cdot w({\cal O}_{h^{k-2-j}}{\cal O}_{h^{j-1}})_{0,d}\cdot e^{du}\right)\no\\
&=&x+\sum_{d=1}^{\infty}\frac{1}{k}w({\cal O}_{h^{k-2-j}}{\cal O}_{h^{j-1}})_{0,d}\cdot e^{dx}.
\ean
The last line is nothing but the assertion of the theorem in the $m=1$ case.

Assuming that the assertion of the theorem holds for $m$, we obtain 
\ban
J_{j,m+1}^{k,k}(x)&=&\int^{x}du\left(1+\sum_{d=1}^{\infty}\frac{d}{k}\cdot w({\cal O}_{h^{k-2-j+m}}{\cal O}_{h^{j-m-1}})_{0,d}\cdot e^{du}\right)\\
&&\times\left(\frac{u^{m}}{m!}+\sum_{n=1}^{m}\frac{u^{m-n}}{(m-n)!}\left(\sum_{f=1}^{\infty}\frac{1}{k}w(\tau_{n-1}({\cal O}_{h^{k-2-j+m-n}}){\cal O}_{h^{j-m}})_{0,f}e^{fu}\right)\right)
\ean
\ban
&=&\int^{x}du\left(\frac{u^{m}}{m!}+\sum_{n=1}^{m}\sum_{f=1}^{\infty}\frac{1}{k}w(\tau_{n-1}({\cal O}_{h^{k-2-j+m-n}}){\cal O}_{h^{j-m}})_{0,f}\frac{u^{m-n}}{(m-n)!}e^{fu}\right)\no\\
&&\times\left(1+\sum_{d=1}^{\infty}\frac{d}{k}\cdot w({\cal O}_{h^{k-2-j+m}}{\cal O}_{h^{j-m-1}})_{0,d}\cdot e^{du}\right)\\
&=&\frac{x^{m+1}}{(m+1)!}+\sum_{n=1}^{m}\sum_{f=1}^{\infty}\frac{1}{k}w(\tau_{n-1}({\cal O}_{h^{k-2-j+m-n}}){\cal O}_{h^{j-m}})_{0,f}\left(\int^{x}du\frac{u^{m-n}}{(m-n)!}e^{fu}\right)\no\\
&&+\sum_{d=1}^{\infty}\frac{d}{k}\cdot w({\cal O}_{h^{k-2-j+m}}{\cal O}_{h^{j-m-1}})_{0,d}\left(\int^{x}du\frac{u^{m}}{m!}e^{du}\right)\no\\
&&+\sum_{n=1}^{m}\sum_{f=1}^{\infty}\sum_{d=1}^{\infty}\frac{1}{k}w(\tau_{n-1}({\cal O}_{h^{k-2-j+m-n}}){\cal O}_{h^{j-m}})_{0,f}\frac{d}{k}\cdot w({\cal O}_{h^{k-2-j+m}}{\cal O}_{h^{j-m-1}})_{0,d}\no\\
&&\times\left(\int^{x}du\frac{u^{m-n}}{(m-n)!}e^{(f+d)u}\right).
\ean 
Therefore, we have
\ban
J_{j,m+1}^{k,k}(x)&=&\frac{x^{m+1}}{(m+1)!}+A+B+C,
\ean
where 
\ban
A&=&\sum_{n=1}^{m}\sum_{d=1}^{\infty}\frac{1}{k}w(\tau_{n-1}({\cal O}_{h^{k-2-j+m-n}}){\cal O}_{h^{j-m}})_{0,d}\left(\int^{x}du\frac{u^{m-n}}{(m-n)!}e^{du}\right)\no\\
&=&\sum_{n=1}^{m}\sum_{d=1}^{\infty}\frac{1}{k}w(\tau_{n-1}({\cal O}_{h^{k-2-j+m-n}}){\cal O}_{h^{j-m}})_{0,d}\left(\sum_{l=0}^{m-n}(-1)^{l}\frac{x^{m-n-l}}{(m-n-l)!}\frac{1}{d^{l+1}}e^{dx}\right),
\ean
\ban
B&=&\sum_{d=1}^{\infty}\frac{d}{k}\cdot w({\cal O}_{h^{k-2-j+m}}{\cal O}_{h^{j-m-1}})_{0,d}\left(\int^{x}du\frac{u^{m}}{m!}e^{du}\right)\no\\
&=&\sum_{d=1}^{\infty}\frac{d}{k}\cdot w({\cal O}_{h^{k-2-j+m}}{\cal O}_{h^{j-m-1}})_{0,d}\left(\sum_{l=0}^{m}(-1)^{l}\frac{x^{m-l}}{(m-l)!}\frac{1}{d^{l+1}}e^{dx}\right),
\ean
\ban
C&=&\sum_{n=1}^{m}\sum_{f=1}^{\infty}\sum_{d=1}^{\infty}\frac{1}{k}w(\tau_{n-1}({\cal O}_{h^{k-2-j+m-n}}){\cal O}_{h^{j-m}})_{0,f}\frac{d}{k}\cdot w({\cal O}_{h^{k-2-j+m}}{\cal O}_{h^{j-m-1}})_{0,d}\no\\
&&\times\left(\int^{x}du\frac{u^{m-n}}{(m-n)!}e^{(f+d)u}\right)\\
&=&\sum_{n=1}^{m}\sum_{d=1}^{\infty}\left(\sum_{f=1}^{d-1}(d-f)\frac{1}{k}w(\tau_{n-1}({\cal O}_{h^{k-2-j+m-n}}){\cal O}_{h^{j-m}})_{0,f}\frac{1}{k}\cdot w({\cal O}_{h^{k-2-j+m}}{\cal O}_{h^{j-m-1}})_{0,d-f}\right)\no\\
&&\times\left(\int^{x}du\frac{u^{m-n}}{(m-n)!}e^{du}\right)\\
&=&\sum_{n=1}^{m}\sum_{d=1}^{\infty}\biggl(d\frac{1}{k}w(\tau_{n}({\cal O}_{h^{k-2-j+m-n}}){\cal O}_{h^{j-m-1}})_{0,d}+\frac{1}{k}w(\tau_{n-1}({\cal O}_{h^{k-1-j+m-n}}){\cal O}_{h^{j-m-1}})_{0,d}
\ean
\ban
&&-\frac{1}{k}w(\tau_{n-1}({\cal O}_{h^{k-2-j+m-n}}){\cal O}_{h^{j-m}})_{0,d}\biggr)\left(\sum_{l=0}^{m-n}(-1)^{l}\frac{x^{m-n-l}}{(m-n-l)!}\frac{1}{d^{l+1}}e^{dx}\right).
\ean
In the above, we used (\ref{int}) and Lemma \ref{wtrl}. We further decompose $C$ into 
three parts:
\ban
C=C1+C2+C3,
\ean
where
\ban
C1&=&\sum_{n=1}^{m}\sum_{d=1}^{\infty}d\frac{1}{k}w(\tau_{n}({\cal O}_{h^{k-2-j+m-n}}){\cal O}_{h^{j-m-1}})_{0,d}\left(\sum_{l=0}^{m-n}(-1)^{l}\frac{x^{m-n-l}}{(m-n-l)!}\frac{1}{d^{l+1}}e^{dx}\right)\\
&=&\sum_{n=1}^{m}\sum_{d=1}^{\infty}\frac{1}{k}w(\tau_{n}({\cal O}_{h^{k-2-j+m-n}}){\cal O}_{h^{j-m-1}})_{0,d}\left(\sum_{l=0}^{m-n}(-1)^{l}\frac{x^{m-n-l}}{(m-n-l)!}\frac{1}{d^{l}}e^{dx}\right),\\
C2&=&\sum_{n=1}^{m}\sum_{d=1}^{\infty}\frac{1}{k}w(\tau_{n-1}({\cal O}_{h^{k-1-j+m-n}}){\cal O}_{h^{j-m-1}})_{0,d}\left(\sum_{l=0}^{m-n}(-1)^{l}\frac{x^{m-n-l}}{(m-n-l)!}\frac{1}{d^{l+1}}e^{dx}\right),\\
C3&=&-\sum_{n=1}^{m}\sum_{d=1}^{\infty}\frac{1}{k}w(\tau_{n-1}({\cal O}_{h^{k-2-j+m-n}}){\cal O}_{h^{j-m}})_{0,d}\left(\sum_{l=0}^{m-n}(-1)^{l}\frac{x^{m-n-l}}{(m-n-l)!}\frac{1}{d^{l+1}}e^{dx}\right).
\ean 
From this decomposition, we can easily see that
\ban
A+C3=0.
\ean 
Next, we rewrite C1 by replacing $n$ with $n-1$:
\ban
C1=\sum_{n=2}^{m+1}\sum_{d=1}^{\infty}\frac{1}{k}w(\tau_{n-1}({\cal O}_{h^{k-1-j+m-n}}){\cal O}_{h^{j-m-1}})_{0,d}\left(\sum_{l=0}^{m+1-n}(-1)^{l}\frac{x^{m+1-n-l}}{(m+1-n-l)!}\frac{1}{d^{l}}e^{dx}\right).
\ean
We further decompose C2 into
\ban
C2=C2A+C2B,
\ean
where
\ban
C2A&=&\sum_{n=2}^{m}\sum_{d=1}^{\infty}\frac{1}{k}w(\tau_{n-1}({\cal O}_{h^{k-1-j+m-n}}){\cal O}_{h^{j-m-1}})_{0,d}\left(\sum_{l=0}^{m-n}(-1)^{l}\frac{x^{m-n-l}}{(m-n-l)!}\frac{1}{d^{l+1}}e^{dx}\right)\\
&=&\sum_{n=2}^{m}\sum_{d=1}^{\infty}\frac{1}{k}w(\tau_{n-1}({\cal O}_{h^{k-1-j+m-n}}){\cal O}_{h^{j-m-1}})_{0,d}\left(\sum_{l=1}^{m+1-n}(-1)^{l-1}\frac{x^{m+1-n-l}}{(m+1-n-l)!}\frac{1}{d^{l}}e^{dx}\right)\\
C2B&=&\sum_{d=1}^{\infty}\frac{1}{k}w({\cal O}_{h^{k-2-j+m}}{\cal O}_{h^{j-m-1}})_{0,d}\left(\sum_{l=0}^{m-1}(-1)^{l}\frac{x^{m-l-1}}{(m-l-1)!}\frac{1}{d^{l+1}}e^{dx}\right)\\
&=&\sum_{d=1}^{\infty}\frac{1}{k}w({\cal O}_{h^{k-2-j+m}}{\cal O}_{h^{j-m-1}})_{0,d}\left(\sum_{l=1}^{m}(-1)^{l-1}\frac{x^{m-l}}{(m-l)!}\frac{1}{d^{l}}e^{dx}\right).
\ean
In the above, we rewrite each term by replacing $l$ with $l-1$. Then we can easily see that
\ban
B+C2B=\sum_{d=1}^{\infty}\frac{1}{k}w({\cal O}_{h^{k-2-j+m}}{\cal O}_{h^{j-m-1}})_{0,d}\frac{x^{m}}{m!}e^{dx}.
\ean
Hence we obtain
\ban
&&J_{j,m+1}^{k,k}(x)\\
&=&\frac{x^{m+1}}{(m+1)!}+C1+C2A+\sum_{d=1}^{\infty}\frac{1}{k}w({\cal O}_{h^{k-2-j+m}}{\cal O}_{h^{j-m-1}})_{0,d}\frac{x^{m}}{m!}e^{dx}\\
&=&\frac{x^{m+1}}{(m+1)!}+\sum_{n=2}^{m+1}\sum_{d=1}^{\infty}\frac{1}{k}w(\tau_{n-1}({\cal O}_{h^{k-1-j+m-n}}){\cal O}_{h^{j-m-1}})_{0,d}\left(\sum_{l=0}^{m+1-n}(-1)^{l}\frac{x^{m+1-n-l}}{(m+1-n-l)!}\frac{1}{d^{l}}e^{dx}\right)\\
&&+\sum_{n=2}^{m}\sum_{d=1}^{\infty}\frac{1}{k}w(\tau_{n-1}({\cal O}_{h^{k-1-j+m-n}}){\cal O}_{h^{j-m-1}})_{0,d}\left(\sum_{l=1}^{m+1-n}(-1)^{l-1}\frac{x^{m+1-n-l}}{(m+1-n-l)!}\frac{1}{d^{l}}e^{dx}\right)\\
&&+\sum_{d=1}^{\infty}\frac{1}{k}w({\cal O}_{h^{k-2-j+m}}{\cal O}_{h^{j-m-1}})_{0,d}\frac{x^{m}}{m!}e^{dx}
\ean
\ban
&=&\frac{x^{m+1}}{(m+1)!}+\sum_{d=1}^{\infty}\frac{1}{k}w(\tau_{m}({\cal O}_{h^{k-2-j}}){\cal O}_{h^{j-m-1}})_{0,d}e^{dx}\\
&&+\sum_{n=2}^{m}\sum_{d=1}^{\infty}\frac{1}{k}w(\tau_{n-1}({\cal O}_{h^{k-1-j+m-n}}){\cal O}_{h^{j-m-1}})_{0,d}\left(\sum_{l=0}^{m+1-n}(-1)^{l}\frac{x^{m+1-n-l}}{(m+1-n-l)!}\frac{1}{d^{l}}e^{dx}\right)\\
&&+\sum_{n=2}^{m}\sum_{d=1}^{\infty}\frac{1}{k}w(\tau_{n-1}({\cal O}_{h^{k-1-j+m-n}}){\cal O}_{h^{j-m-1}})_{0,d}\left(\sum_{l=1}^{m+1-n}(-1)^{l-1}\frac{x^{m+1-n-l}}{(m+1-n-l)!}\frac{1}{d^{l}}e^{dx}\right)\\
&&+\sum_{d=1}^{\infty}\frac{1}{k}w({\cal O}_{h^{k-2-j+m}}{\cal O}_{h^{j-m-1}})_{0,d}\frac{x^{m}}{m!}e^{dx}
\ean
\ban
&=&\frac{x^{m+1}}{(m+1)!}+\sum_{d=1}^{\infty}\frac{1}{k}w(\tau_{m}({\cal O}_{h^{k-2-j}}){\cal O}_{h^{j-m-1}})_{0,d}e^{dx}\\
&&+\sum_{n=2}^{m}\sum_{d=1}^{\infty}\frac{1}{k}w(\tau_{n-1}({\cal O}_{h^{k-1-j+m-n}}){\cal O}_{h^{j-m-1}})_{0,d}\frac{x^{m+1-n}}{(m+1-n)!}e^{dx}\\
&&+\sum_{d=1}^{\infty}\frac{1}{k}w({\cal O}_{h^{k-2-j+m}}{\cal O}_{h^{j-m-1}})_{0,d}\frac{x^{m}}{m!}e^{dx}
\ean
\ban
&=&\frac{x^{m+1}}{(m+1)!}+\sum_{n=1}^{m+1}\sum_{d=1}^{\infty}\frac{1}{k}w(\tau_{n-1}({\cal O}_{h^{k-2-j+m+1-n}}){\cal O}_{h^{j-m-1}})_{0,d}\frac{x^{m+1-n}}{(m+1-n)!}e^{dx}.
\ean
This completes the induction. $\Box$

\begin{lem}
\ba
&&\sum_{f=1}^{d-1}(d-f)\frac{1}{k}\langle\tau_{n-1}({\cal O}_{h^{k-2-j+m-n}}){\cal O}_{h^{j-m}}\rangle_{0,f}\cdot
\frac{1}{k}\langle{\cal O}_{h^{k-2-j+m}}{\cal O}_{h^{j-m-1}}\rangle_{0,d-f}\no\\
&=&d\cdot\frac{1}{k}\langle\tau_{n}({\cal O}_{h^{k-2-j+m-n}}){\cal O}_{h^{j-m-1}}\rangle_{0,d}+
\frac{1}{k}\langle\tau_{n-1}({\cal O}_{h^{k-1-j+m-n}}){\cal O}_{h^{j-m-1}}\rangle_{0,d}\no\\
&&-\frac{1}{k}\langle\tau_{n-1}({\cal O}_{h^{k-2-j+m-n}}){\cal O}_{h^{j-m}}\rangle_{0,d}.
\label{trr}
\ea
\label{trl}
\end{lem}
\noindent
{\bf Proof)} From the topological recursion relation, we obtain
\ba
&&\sum_{f=0}^{d}\langle\tau_{n-1}({\cal O}_{h^{k-2-j+m-n}}){\cal O}_{h^{j-m}}\rangle_{0,f}\cdot
\frac{1}{k}\langle{\cal O}_{h^{k-2-j+m}}{\cal O}_{h^{j-m-1}}{\cal O}_{h}\rangle_{0,d-f}\no\\
&=&\langle\tau_{n}({\cal O}_{h^{k-2-j+m-n}}){\cal O}_{h^{j-m-1}}{\cal O}_{h}\rangle_{0,d}
\label{lemeq}
\ea
A degree $0$ gravitational two-point GW invariant always vanishes, $\langle{\cal O}_{h^{k-2-j+m}}{\cal O}_{h^{j-m-1}}{\cal O}_{h}\rangle_{0,0}=k$, and $\langle{\cal O}_{h^{k-2-j+m}}{\cal O}_{h^{j-m-1}}{\cal O}_{h}\rangle_{0,d-f}=(d-f)\langle{\cal O}_{h^{k-2-j+m}}{\cal O}_{h^{j-m-1}}\rangle_{0,d-f}$ by the divisor axiom. Hence, the l.h.s. of (\ref{lemeq}) is rewritten as follows:
\ba
&&\sum_{f=0}^{d-1}(d-f)\langle\tau_{n-1}({\cal O}_{h^{k-2-j+m-n}}){\cal O}_{h^{j-m}}\rangle_{0,f}\cdot
\frac{1}{k}\langle{\cal O}_{h^{k-2-j+m}}{\cal O}_{h^{j-m-1}}\rangle_{0,d-f}\no\\
&&+\langle\tau_{n-1}({\cal O}_{h^{k-2-j+m-n}}){\cal O}_{h^{j-m}}\rangle_{0,d}.
\label{lhs}
\ea
On the other hand, by applying Hori's equation, the r.h.s. becomes 
\ba
d\langle\tau_{n}({\cal O}_{h^{k-2-j+m-n}}){\cal O}_{h^{j-m-1}}\rangle_{0,d}+\langle\tau_{n-1}({\cal O}_{h^{k-1-j+m-n}}){\cal O}_{h^{j-m-1}}\rangle_{0,d}
\label{rhs}
\ea
By equating (\ref{lhs}) and (\ref{rhs}), we obtain the desired equality. \quad  $\Box$
\begin{defi}
\ban
w({\cal O}_{h^{k-2-j}}{\cal O}_{h^{j-1}})_{0}(x)&:=&kx+\sum_{d=1}^{\infty}w({\cal O}_{h^{k-2-j}}{\cal O}_{h^{j-1}})_{0,d}e^{dx},\\
\langle{\cal O}_{h^{k-2-j}}{\cal O}_{h^{j-1}}\rangle_{0}(t)&:=&kt+\sum_{d=1}^{\infty}\langle{\cal O}_{h^{k-2-j}}{\cal O}_{h^{j-1}}\rangle_{0,d}e^{dt},\\
\ean
\end{defi}
With this definition, we recall the mirror theorem for the virtual structure constants \cite{Jin6}.
\begin{theorem} {\bf \cite{Jin6}}
The following equality holds:
\ban
\langle{\cal O}_{h^{k-2-j}}{\cal O}_{h^{j-1}}\rangle_{0}(t(x))&=&w({\cal O}_{h^{k-2-j}}{\cal O}_{h^{j-1}})_{0}(x)\no\\
(\Longleftrightarrow\;\;\langle{\cal O}_{h^{k-2-j}}{\cal O}_{h^{j-1}}\rangle_{0}(t)&=&w({\cal O}_{h^{k-2-j}}{\cal O}_{h^{j-1}})_{0}(x(t))\;\;),
\ean
where $t=t(x)$ ($x=x(t)$) is the mirror map (and its inverse) given in Definition \ref{mmap}.
\label{vmth}
\end{theorem}
\vspace{1em}
\noindent
{\bf Proof of Theorem \ref{mainA})} We also prove the theorem by induction on $m$. In the $m=1$ case, 
the assertion of the theorem corresponds to
\ban
\frac{1}{k}w({\cal O}_{h^{k-2-j}}{\cal O}_{h^{j-1}})_{0}(x(t))=\frac{1}{k}\langle{\cal O}_{h^{k-2-j}}{\cal O}_{h^{j-1}}\rangle_{0}(t),
\ean 
which follows from Theorem \ref{vmth}.

Next, we assume that the theorem holds for $m$. By applying integration by substitution to the definition of $J_{j,m+1}^{k,k}(x)$, we obtain
\ban
&&J_{j,m+1}^{k,k}(x(t))\\&=&\int^{t}\tilde{L}_{j-m}^{k,k}(x(t))\left(\frac{t^{m}}{m!}+\sum_{n=1}^{m}\frac{t^{m-n}}{(m-n)!}\left(\sum_{d=1}^{\infty}\frac{1}{k}\langle\tau_{n-1}({\cal O}_{h^{k-2-j+m-n}}){\cal O}_{h^{j-m}}\rangle_{0,d}e^{dt}\right)\right)\frac{dx(t)}{dt}dt.\\
\ean
Notice here that $\displaystyle{\tilde{L}_{j-m}^{k,k}(x(t))=\frac{1}{k}\frac{dw({\cal O}_{h^{k-2-j+m}}{\cal O}_{h^{j-m-1}})_{0}}{dx}(x(t))}$. Hence we obtain
\ban
&&J_{j,m+1}^{k,k}(x(t))\\
&=&\int^{t}\left(\frac{t^{m}}{m!}+\sum_{n=1}^{m}\frac{t^{m-n}}{(m-n)!}\left(\sum_{f=1}^{\infty}\frac{1}{k}\langle\tau_{n-1}({\cal O}_{h^{k-2-j+m-n}}){\cal O}_{h^{j-m}}\rangle_{0,f}e^{ft}\right)\right)
\\
&&\times\frac{1}{k}\frac{dw({\cal O}_{h^{k-2-j+m}}{\cal O}_{h^{j-m-1}})_{0}}{dx}(x(t))\frac{dx(t)}{dt}dt\\
&=&\int^{t}\left(\frac{t^{m}}{m!}+\sum_{n=1}^{m}\frac{t^{m-n}}{(m-n)!}\left(\sum_{f=1}^{\infty}\frac{1}{k}\langle\tau_{n-1}({\cal O}_{h^{k-2-j+m-n}}){\cal O}_{h^{j-m}}\rangle_{0,f}e^{ft}\right)\right)\\
&&\times\frac{1}{k}\left(\frac{d}{dt}w({\cal O}_{h^{k-2-j+m}}{\cal O}_{h^{j-m-1}})_{0}(x(t))\right)dt.\\
\ean
Then using Theorem \ref{vmth} again, we can rewrite $J_{j,m+1}^{k,k}(x(t))$ into the form:
\ban
&&J_{j,m+1}^{k,k}(x(t))\\
&=&\int^{t}\left(\frac{t^{m}}{m!}+\sum_{n=1}^{m}\frac{t^{m-n}}{(m-n)!}\left(\sum_{f=1}^{\infty}\frac{1}{k}\langle\tau_{n-1}({\cal O}_{h^{k-2-j+m-n}}){\cal O}_{h^{j-m}}\rangle_{0,f}e^{ft}\right)\right)\\
&&\times\frac{1}{k}\left(\frac{d}{dt}\langle{\cal O}_{h^{k-2-j+m}}{\cal O}_{h^{j-m-1}}\rangle_{0}(t)\right)dt\\
&=&\int^{t}\left(\frac{t^{m}}{m!}+\sum_{n=1}^{m}\frac{t^{m-n}}{(m-n)!}\left(\sum_{f=1}^{\infty}\frac{1}{k}\langle\tau_{n-1}({\cal O}_{h^{k-2-j+m-n}}){\cal O}_{h^{j-m}}\rangle_{0,f}e^{ft}\right)\right)\\
&&\times\left(1+\sum_{d=1}^{\infty}\frac{d}{k}\langle{\cal O}_{h^{k-2-j+m}}{\cal O}_{h^{j-m-1}}\rangle_{0,d}e^{dt}\right)dt.
\ean
Using Lemma \ref{trl}, the remaining computation proceeds in exactly the same way as in the proof of Theorem
\ref{mainB}, and we obtain
\ban
J_{j,m+1}^{k,k}(x(t))=\frac{t^{m+1}}{(m+1)!}+\sum_{n=1}^{m+1}\sum_{d=1}^{\infty}\frac{1}{k}\langle\tau_{n-1}({\cal O}_{h^{k-2-j+m+1-n}}){\cal O}_{h^{j-m-1}}\rangle_{0,d}\frac{t^{m+1-n}}{(m+1-n)!}e^{dt},
\ean
which completes the induction. \quad $\Box$
\section{Application of the Main Theorems}
\subsection{Connection with Givental's $J$-function}
First, we discuss the relation between Theorem \ref{mainB} and the generalized hypergeometric functions 
used in the mirror computation of the Calabi-Yau hypersurface $M_{k}^{k}$. Again using Lemma \ref{wtrl}, 
we can derive the following equality:
\ba
&&J_{j,j}^{k,k}(x)\tilde{L}_{0}^{k,k}(x)=\frac{x^{j}}{j!}\left(1+\sum_{d=1}^{\infty}\frac{d}{k}w({\cal O}_{h^{k-2}}{\cal O}_{h^{-1}})_{0,d}e^{dx} \right)\no\\
&&+\sum_{n=1}^{j}\frac{x^{j-n}}{(j-n)!}\sum_{d=1}^{\infty}\frac{1}{k}
\biggl(dw(\tau_{n}({\cal O}_{h^{k-2-n}}){\cal O}_{h^{-1}})_{0,d}+w(\tau_{n-1}({\cal O}_{h^{k-1-n}}){\cal O}_{h^{-1}})_{0,d}\biggr)e^{dx}.\no\\
\label{JL}
\ea
At this stage, let us recall the main theorem of \cite{JM}.
\begin{theorem}{\bf \cite{JM}}
\ba
&&\frac{d}{k}w(\tau_{n}({\cal O}_{h^{k-2-n}}){\cal O}_{h^{-1}})_{0,d}+
\frac{1}{k}w(\tau_{n-1}({\cal O}_{h^{k-1-n}}){\cal O}_{h^{-1}})_{0,d}=\frac{1}{n!}\frac{\d^{n}}{\d \epsilon^{n}}\left(\frac{\prod_{r=1}^{kd}(r+k\epsilon)}{\prod_{r=1}^{d}(r+\epsilon)^k}\middle) \right|_{\epsilon=0}.\;(n\geq 1),\no\\
&&\frac{d}{k}w({\cal O}_{h^{k-2}}{\cal O}_{h^{-1}})_{0,d}=\frac{(kd)!}{(d!)^{k}}.
\label{mainth}
\ea
\end{theorem}
Then (\ref{JL}) turns into
\ba
&&J_{j,j}^{k,k}(x)\tilde{L}_{0}^{k,k}(x)\no\\
&=&\frac{x^{j}}{j!}\left(\sum_{d=0}^{\infty}\frac{(kd)!}{(d!)^{k}}e^{dx} \right)
+\sum_{n=1}^{j}\frac{x^{j-n}}{(j-n)!n!}\sum_{d=1}^{\infty}\frac{\d^{n}}{\d \epsilon^{n}}\left(\frac{\prod_{r=1}^{kd}(r+k\epsilon)}{\prod_{r=1}^{d}(r+\epsilon)^k}\middle) \right|_{\epsilon=0}e^{dx}\no\\
&=&\frac{1}{j!}\sum_{n=0}^{j}{j \choose n}\frac{\d^{n}}{\d \epsilon^{n}}
\left(\frac{\prod_{r=1}^{kd}(r+k\epsilon)}{\prod_{r=1}^{d}(r+\epsilon)^k}\middle) \right|_{\epsilon=0}e^{dx}.
\ea
The expression in the last line:
\ba
W^{k,k}_{(j)}(x):=\sum_{n=0}^{j}{j \choose n}\frac{\d^{n}}{\d \epsilon^{n}}\left(\frac{\prod_{r=1}^{kd}(r+k\epsilon)}{\prod_{r=1}^{d}(r+\epsilon)^k}\middle) \right|_{\epsilon=0}e^{dx}\quad(j=0,1,\cdots,k-1),
\ea
is the $j$-th log-solution of the differential equation:
\ba
\left(\frac{d}{dx}\biggl((\frac{d}{dx})^{k-1}-k\cdot e^{x}\cdot (k\frac{d}{dx}+(k-1))(k\frac{d}{dx}+(k-2))\cdots(k\frac{d}{dx}+1)\biggr)\right)W(x)=0.\no\\
\label{eq2}
\ea
By noting the equality $\tilde{L}_{0}^{k,k}(x)=W^{k,k}_{(0)}(x)$, we obtain
\begin{cor}
\ban
J_{j,j}^{k,k}(x)W^{k,k}_{(0)}(x)&=&\frac{1}{j!}W^{k,k}_{(j)}(x)\\
\Longleftrightarrow J_{j,j}^{k,k}(x)&=&\frac{1}{j!}\frac{W^{k,k}_{(j)}(x)}{W^{k,k}_{(0)}(x)} \quad 
(0\leq j\leq  k-1),
\ean
where 
\ban
W^{k,k}_{(j)}(x):=\sum_{n=0}^{j}{j \choose n}\frac{\d^{n}}{\d \epsilon^{n}}\left(\frac{\prod_{r=1}^{kd}(r+k\epsilon)}{\prod_{r=1}^{d}(r+\epsilon)^k}\middle) \right|_{\epsilon=0}e^{dx}\quad(j=0,1,\cdots,k-1)       ,
\ean
is the $j$-th log-solution of the differential equation:
\ban
\left(\frac{d}{dx}\biggl((\frac{d}{dx})^{k-1}-k\cdot e^{x}\cdot (k\frac{d}{dx}+(k-1))(k\frac{d}{dx}+(k-2))\cdots(k\frac{d}{dx}+1)\biggr)\right)W(x)=0.
\ean
\label{Jcor}
\end{cor}
Combination of the above corollary with Theorem \ref{mainA} results in  
\begin{cor}
\ban
J_{j,j}^{k,k}(x(t))&=&\frac{1}{j!}\frac{W^{k,k}_{(j)}(x(t))}{W^{k,k}_{(0)}(x(t))}\\
                        &=&\frac{t^{j}}{j!}+\sum_{n=1}^{j}\frac{t^{j-n}}{(j-n)!}\left(\sum_{d=1}^{\infty}\frac{1}{k}\langle\tau_{n-1}({\cal O}_{h^{k-2-n}}){\cal O}_{h^{0}}\rangle_{0,d}e^{dt}\right),
\ean
where $x=x(t)$ is the inverse of the mirror map in Definition \ref{mmap}.
\end{cor}
{\bf This corollary corresponds to the statement of Givental's version of the Mirror Theorem \cite{Giv} in terms of the $J$-function.} Therefore, we have proved Givental's version of the Mirror Theorem from the point of view 
of virtual structure constants. In our context, Givental's $I$-function is given by
\ba
I(x):=\sum_{j=0}^{k-2}\frac{1}{j!}W^{k,k}_{(j)}(x)P^{j},
\ea
and the $J$-function corresponds to
\ba
J(t)&:=&\frac{I(x(t))}{W^{k,k}_{(0)}(x(t))}=\sum_{j=0}^{k-2}\frac{1}{j!}\frac{W^{k,k}_{(j)}(x(t))}{W^{k,k}_{(0)}(x(t))}P^{j}\no\\
&=&\sum_{j=0}^{k-2}J_{j,j}^{k,k}(x(t))P^{j}.
\ea
By tracing back our construction, we can represent $\tilde{L}_{j}^{k,k}(x)$ and $J_{j,m}^{k,k}(x)$ in terms of $W^{k,k}_{(j)}(x)\;\;(j=0,\cdots,k-1)$.
\begin{cor}
$\tilde{L}_{j}^{k,k}(x)$ is inductively determined from $W^{k,k}_{(j)}(x)\;\;(j=0,\cdots,k-1)$ by 
the initial condition:
\ban
\tilde{L}_{0}^{k,k}(x)=W^{k,k}_{(0)}(x),
\ean 
and the equality:
\ban
\tilde{L}_{j}^{k,k}(x)=\frac{1}{j!}\frac{d}{dx}\left(\frac{1}{\tilde{L}_{j-1}^{k,k}(x)}\frac{d}{dx}\left(
\frac{1}{\tilde{L}_{j-2}^{k,k}(x)}\frac{d}{dx}\left(\cdots\frac{1}{\tilde{L}_{1}^{k,k}(x)}\frac{d}{dx}\left(
\frac{W_{(j)}^{k,k}(x)}{\tilde{L}_{0}^{k,k}(x)}\right)\cdots\right)\right)\right).
\ean
Then $J_{j,m}^{k,k}(x)$ is determined by the formula:
\ban
J_{j,m}^{k,k}(x)=\frac{1}{j!}\left(\frac{1}{\tilde{L}_{j-m}^{k,k}(x)}\frac{d}{dx}\left(
\frac{1}{\tilde{L}_{j-m-1}^{k,k}(x)}\frac{d}{dx}\left(\cdots\frac{1}{\tilde{L}_{1}^{k,k}(x)}\frac{d}{dx}\left(
\frac{W_{(j)}^{k,k}(x)}{\tilde{L}_{0}^{k,k}(x)}\right)\cdots\right)\right)\right).
\ean
\label{fin}
\end{cor}
By combining this corollary with Theorem \ref{mainA}, we can determine all the \\$\langle\tau_{n-1}({\cal O}_{h^{k-2-j+m-n}}){\cal O}_{h^{j-m}}\rangle_{0,d}$'s of $M_{k}^{k}$ from $W^{k,k}_{(j)}(x)\;\;(j=0,\cdots,k-1)$, which corresponds to {\bf the classical mirror computation of the Calabi-Yau hypersurface $M_{k}^{k}$} \cite{bara, gmp, JN}.

\begin{Rem}
The construction of this paper is motivated by Theorem 4.2 in \cite{Jin1}, where the following equality
is given:
\ban
&&\left(\frac{d}{dx}\left(\frac{1}{\tilde{L}_{k-1}^{k,k}(x)}\frac{d}{dx}\left(
\frac{1}{\tilde{L}_{k-2}^{k,k}(x)}\frac{d}{dx}\left(\cdots\frac{1}{\tilde{L}_{1}^{k,k}(x)}\frac{d}{dx}\left(
\frac{W(x)}{\tilde{L}_{0}^{k,k}(x)}\right)\cdots\right)\right)\right)\right)\\
&=&\left(\frac{d}{dx}\biggl((\frac{d}{dx})^{k-1}-k\cdot e^{x}\cdot (k\frac{d}{dx}+(k-1))(k\frac{d}{dx}+(k-2))\cdots(k\frac{d}{dx}+1)\biggr)\right)W(x).
\ean
\end{Rem}

\subsection{Examples and Investigation of ``Ghost Gromov-Witten Invariants''}
\begin{ex} {\bf Quintic Hypersurface $M_{5}^{5}$}
In this case, the non-trivial $J$'s are
\ban
&&J_{1,1}^{5,5},J_{2,2}^{5,5}, J_{2,1}^{5,5},J_{3,3}^{5,5},J_{3,2}^{5,5},J_{3,1}^{5,5},J_{4,4,}^{5,5},J_{4,3}^{5,5},J_{4,2}^{5,5},J_{4,1}^{5,5}.\\
\ean 
Theorem \ref{mainA} tells us that they have the following structure up to $j=3$:
\ban
5J_{1,1}^{5,5}(x(t))&=&5t+\sum_{d=1}^{\infty}\langle{\cal O}_{h^2}{\cal O}_{h^{0}}\rangle_{0,d}e^{dt}=5t,\\
5J_{2,2}^{5,5}(x(t))&=&5\frac{t^2}{2}+t\left(\sum_{d=1}^{\infty}\langle{\cal O}_{h}{\cal O}_{h^{0}}\rangle_{0,d}e^{dt}\right)+\sum_{d=1}^{\infty}\langle{\tau_{1}(\cal O}_{h}){\cal O}_{h^{0}}\rangle_{0,d}e^{dt}=5\frac{t^2}{2}+\sum_{d=1}^{\infty}\langle{\cal O}_{h}\rangle_{0,d}e^{dt},
\ean
\ban
5J_{2,1}^{5,5}(x(t))&=&5t+\sum_{d=1}^{\infty}\langle{\cal O}_{h}{\cal O}_{h}\rangle_{0,d}e^{dt},
\ean
\ban
5J_{3,3}^{5,5}(x(t))&=&5\frac{t^3}{6}+\frac{t^{2}}{2}\left(\sum_{d=1}^{\infty}\langle{\cal O}_{h^2}{\cal O}_{h^{0}}\rangle_{0,d}e^{dt}\right)+t\left(\sum_{d=1}^{\infty}\langle\tau_{1}({\cal O}_{h}){\cal O}_{h^{0}}\rangle_{0,d}e^{dt}\right)\\
&&+\sum_{d=1}^{\infty}\langle\tau_{2}({\cal O}_{h^{0}}){\cal O}_{h^{0}}\rangle_{0,d}e^{dt}\\
&=&5\frac{t^3}{6}+t\left(\sum_{d=1}^{\infty}\langle{\cal O}_{h}\rangle_{0,d}e^{dt}\right)-2\sum_{d=1}^{\infty}\langle*\rangle_{0,d}e^{dt},
\ean
\ban
5J_{3,2}^{5,5}(x(t))&=&5\frac{t^2}{2}+t\left(\sum_{d=1}^{\infty}\langle{\cal O}_{h}{\cal O}_{h}\rangle_{0,d}e^{dt}\right)+\sum_{d=1}^{\infty}\langle\tau_{1}({\cal O}_{h^{0}}){\cal O}_{h}\rangle_{0,d}e^{dt}\\
&=&5\frac{t^2}{2}+t\left(\sum_{d=1}^{\infty}\langle{\cal O}_{h}{\cal O}_{h}\rangle_{0,d}e^{dt}\right)-\sum_{d=1}^{\infty}\langle{\cal O}_{h}\rangle_{0,d}e^{dt},\\
5J_{3,1}^{5,5}(x(t))&=&5t+\sum_{d=1}^{\infty}\langle{\cal O}_{h^{0}}{\cal O}_{h^2}\rangle_{0,d}e^{dt}=5t.
\ean
In the above, we used the properties of gravitational Gromov-Witten invariants (\ref{punc}), (\ref{dila}),
and (\ref{hori}). 

If we introduce the ``prepotential'', i.e., the generating function of the Gromov-Witten invariants:
\ban
F(t):=5\frac{t^{3}}{6}+\sum_{d=1}^{\infty}\langle*\rangle_{0,d}e^{dt},
\ean
we obtain
\ban
5J_{2,2}^{5,5}(x(t))=\frac{d}{dt}F(t),\;\;5J_{3,3}^{5,5}(x(t))=t\frac{d}{dt}F(t)-2F(t),
\ean
which corresponds to the well-known classical result.

Next, we enter an exotic domain, namely the $j=4$ case. In this case, Theorem \ref{mainA} tells us
\ban
5J_{4,4}^{5,5}(x(t))&=&5\frac{t^4}{24}+\frac{t^{3}}{6}\left(\sum_{d=1}^{\infty}\langle{\cal O}_{h^2}{\cal O}_{h^{0}}\rangle_{0,d}e^{dt}\right)+\frac{t^2}{2}\left(\sum_{d=1}^{\infty}\langle\tau_{1}({\cal O}_{h}){\cal O}_{h^{0}}\rangle_{0,d}e^{dt}\right)\\
&&+t\left(\sum_{d=1}^{\infty}\langle\tau_{2}({\cal O}_{h^{0}}){\cal O}_{h^{0}}\rangle_{0,d}e^{dt}\right)
+\sum_{d=1}^{\infty}\langle\tau_{3}({\cal O}_{h^{-1}}){\cal O}_{h^{0}}\rangle_{0,d}e^{dt}\\
&=&5\frac{t^4}{24}+\frac{t^2}{2}\left(\sum_{d=1}^{\infty}\langle{\cal O}_{h}\rangle_{0,d}e^{dt}\right)-2t\left(\sum_{d=1}^{\infty}\langle*\rangle_{0,d}e^{dt}\right)
+\sum_{d=1}^{\infty}\langle\tau_{3}({\cal O}_{h^{-1}}){\cal O}_{h^{0}}\rangle_{0,d}e^{dt},
\\
5J_{4,3}^{5,5}(x(t))&=&5\frac{t^3}{6}+\frac{t^2}{2}\left(\sum_{d=1}^{\infty}\langle{\cal O}_{h}{\cal O}_{h}\rangle_{0,d}e^{dt}\right)+t\left(\sum_{d=1}^{\infty}\langle\tau_{1}({\cal O}_{h^{0}}){\cal O}_{h}\rangle_{0,d}e^{dt}\right)\\
&=&5\frac{t^3}{6}+\frac{t^2}{2}\left(\sum_{d=1}^{\infty}\langle{\cal O}_{h}{\cal O}_{h}\rangle_{0,d}e^{dt}\right)-t\left(\sum_{d=1}^{\infty}\langle{\cal O}_{h}\rangle_{0,d}e^{dt}\right)
\ean
\ban
&&+\sum_{d=1}^{\infty}\langle\tau_{2}({\cal O}_{h^{-1}}){\cal O}_{h}\rangle_{0,d}e^{dt},\\
5J_{4,2}^{5,5}(x(t))&=&5\frac{t^2}{2}+t\left(\sum_{d=1}^{\infty}\langle{\cal O}_{h^{0}}{\cal O}_{h^2}\rangle_{0,d}e^{dt}\right)+\sum_{d=1}^{\infty}\langle\tau_{1}({\cal O}_{h^{-1}}){\cal O}_{h^2}\rangle_{0,d}e^{dt}\\
&=&5\frac{t^2}{2}+\sum_{d=1}^{\infty}\langle\tau_{1}({\cal O}_{h^{-1}}){\cal O}_{h^2}\rangle_{0,d}e^{dt},\\
5J_{4,1}^{5,5}(x(t))&=&5t+\sum_{d=1}^{\infty}\langle{\cal O}_{h^{-1}}{\cal O}_{h^3}\rangle_{0,d}e^{dt}.
\ean
In the above, we again used the properties of gravitational Gromov-Witten invariants (\ref{punc}), (\ref{dila}),
and (\ref{hori}). 

The results of explicit computation using Corollary \ref{fin} are given as follows:
\ban
5J_{4,4}^{5,5}(x(t))&=&\frac{5 t^{4}}{24}+\frac{t^2}{2}\left(2875e^t +\frac{4876875}{4} e^{2t}+\frac{8564575000}{9} e^{3t}+\frac{15517926796875}{16} e^{4t}+\cdots\right)\\
&& +t\left(-5750 e^t -\frac{4876875}{4} e^{2t}-\frac{17129150000}{27} e^{3t}-\frac{15517926796875}{32} e^{4t}-\cdots\right) \\
&& +5375 e^t +\frac{5549375 }{16}e^{2t}-\frac{13200321875 }{162}e^{3t}-\frac{474193099015625}{2304}e^{4t}-\cdots,
\ean
\ban
5J_{4,3}^{5,5}(x(t))&=&\frac{5 t^{3}}{6}+ \frac{t^{2}}{2}\left(2875 e^t +\frac{4876875}{2} e^{2t}+\frac{8564575000}{3} e^{3t}+\frac{15517926796875}{4} e^{4t}*\cdots\right)\\
&&+t\left(-2875 e^t -\frac{4876875}{4} e^{2t}-\frac{8564575000}{9} e^{3t}-\frac{15517926796875}{16} e^{4t}-\cdots\right) \\
&& -375 e^t -\frac{4204375}{8}e^{2t}-\frac{15819540625}{18}e^{3t}-\frac{753515781359375}{576}e^{4t}-\cdots,\\
5J_{4,2}^{5,5}(x(t))&=&\frac{5}{2} t^{2}-3250 e^t -\frac{803125}{2} e^{2t}-\frac{1686081250}{9} e^{3t}-\frac{1089994203125}{8} e^{4t}-\cdots,\\
5J_{4,1}^{5,5}(x(t))&=&5t -3250 e^t -803125 e^{2t}-\frac{1686081250}{3} e^{3t}-\frac{1089994203125}{2} e^{4t}-\cdots.
\ean
Let us now look into the exotic parts:
\ba
\sum_{d=1}^{\infty}\langle\tau_{3}({\cal O}_{h^{-1}}){\cal O}_{h^{0}}\rangle_{0,d}e^{dt}
&=&5375 e^t +\frac{5549375 }{16}e^{2t}-\frac{13200321875 }{162}e^{3t}-\frac{474193099015625}{2304}e^{4t}-\cdots,\no\\
\sum_{d=1}^{\infty}\langle\tau_{2}({\cal O}_{h^{-1}}){\cal O}_{h}\rangle_{0,d}e^{dt}
&=& -375 e^t -\frac{4204375}{8}e^{2t}-\frac{15819540625}{18}e^{3t}-\frac{753515781359375}{576}e^{4t}-\cdots,\no\\
\sum_{d=1}^{\infty}\langle\tau_{1}({\cal O}_{h^{-1}}){\cal O}_{h^2}\rangle_{0,d}e^{dt}
&=&-3250 e^t -\frac{803125}{2} e^{2t}-\frac{1686081250}{9} e^{3t}-\frac{1089994203125}{8} e^{4t}-\cdots,\no\\
\sum_{d=1}^{\infty}\langle{\cal O}_{h^{-1}}{\cal O}_{h^3}\rangle_{0,d}e^{dt}
&=&-3250 e^t -803125 e^{2t}-\frac{1686081250}{3} e^{3t}-\frac{1089994203125}{2} e^{4t}-\cdots.\no\\
\label{gdata}
\ea
These are computed from our formal setup, where we use an unnatural 
element $h^{-1}$, i.e., the ``inverse of 
the K\"{a}hler form''. Therefore, we call these rational numbers {\bf ``ghost Gromov-Witten invariants''}. In what follows, we discuss whether they have any geometrical meaning. 
For simplicity, we consider the degree $1$ gravitational ghost Gromov-Witten invariants:
\ban
\langle\tau_{1}({\cal O}_{h^{-1}}){\cal O}_{h^2}\rangle_{0,1}=-3250,\;\;\langle\tau_{2}({\cal O}_{h^{-1}}){\cal O}_{h}\rangle_{0,1}=-375,\;\;\langle\tau_{3}({\cal O}_{h^{-1}}){\cal O}_{h^{0}}\rangle_{0,1}=5375.
\ean
We apply the defining properties of gravitational Gromov-Witten invariants formally to these ghost invariants.
For $\langle\tau_{1}({\cal O}_{h^{-1}}){\cal O}_{h^2}\rangle_{0,1}$, by applying Hori's equation, we obtain
\ban
\langle\tau_{1}({\cal O}_{h^{-1}}){\cal O}_{h^2}{\cal O}_{h}\rangle_{0,1}
&=&\langle\tau_{1}({\cal O}_{h^{-1}}){\cal O}_{h^2}\rangle_{0,1}+\langle{\cal O}_{h^{0
}}{\cal O}_{h^2}\rangle_{0,1}\\
&=&\langle\tau_{1}({\cal O}_{h^{-1}}){\cal O}_{h^2}\rangle_{0,1}.
\ean  
On the other hand, the topological recursion relation gives
\ban
\langle\tau_{1}({\cal O}_{h^{-1}}){\cal O}_{h^2}{\cal O}_{h}\rangle_{0,1}
&=&\frac{1}{5}\langle{\cal O}_{h^{-1}}{\cal O}_{h^3}\rangle_{0,1}\langle{\cal O}_{h^{0}}{\cal O}_{h^2}{\cal O}_{h}\rangle_{0,0}+\frac{1}{5}\langle{\cal O}_{h^{-1}}{\cal O}_{h^3}\rangle_{0,0}\langle
{\cal O}_{h^{0}}{\cal O}_{h^2}{\cal O}_{h}\rangle_{0,1}\\
&=&\langle{\cal O}_{h^{-1}}{\cal O}_{h^3}\rangle_{0,1}.
\ean
Hence we obtain $\langle\tau_{1}({\cal O}_{h^{-1}}){\cal O}_{h^2}\rangle_{0,1}=\langle{\cal O}_{h^{-1}}{\cal O}_{h^3}\rangle_{0,1}$, which agrees with the numerical results in (\ref{gdata}). Next, again from 
Hori's equation, we have
\ban
\langle\tau_{2}({\cal O}_{h^{-1}}){\cal O}_{h}{\cal O}_{h}\rangle_{0,1}
&=&\langle\tau_{2}({\cal O}_{h^{-1}}){\cal O}_{h}\rangle_{0,1}+\langle\tau_{1}({\cal O}_{h^{0
}}){\cal O}_{h}\rangle_{0,1}\\
&=&\langle\tau_{2}({\cal O}_{h^{-1}}){\cal O}_{h}\rangle_{0,1}-\langle{\cal O}_{h}\rangle_{1,0}.
\ean  
In the last line, we used the dilaton equation. We also have the relation from the topological recursion relation:
\ban
\langle\tau_{2}({\cal O}_{h^{-1}}){\cal O}_{h}{\cal O}_{h}\rangle_{0,1}
&=&\frac{1}{5}\langle\tau_{1}({\cal O}_{h^{-1}}){\cal O}_{h^2}\rangle_{0,1}\langle{\cal O}_{h}{\cal O}_{h}{\cal O}_{h}\rangle_{0,0}+\frac{1}{5}\langle\tau_{1}({\cal O}_{h^{-1}}){\cal O}_{h^2}\rangle_{0,0}\langle
{\cal O}_{h}{\cal O}_{h}{\cal O}_{h}\rangle_{0,1}\\
&=&\langle\tau_{1}({\cal O}_{h^{-1}}){\cal O}_{h^2}\rangle_{0,1}.
\ean
Hence we obtain
\ban
\langle\tau_{2}({\cal O}_{h^{-1}}){\cal O}_{h}\rangle_{0,1}=\langle\tau_{1}({\cal O}_{h^{-1}}){\cal O}_{h^2}\rangle_{0,1}+\langle{\cal O}_{h}\rangle_{1,0}=-3250+2875=-375,
\ean
which agrees with the numerical results in (\ref{gdata}). In the same manner, we can derive
\ban
\langle\tau_{3}({\cal O}_{h^{-1}}){\cal O}_{h^{0}}\rangle_{0,1}=\langle\tau_{2}({\cal O}_{h^{-1}}){\cal O}_{h}\rangle_{0,1}+2\langle*\rangle_{1,0}=-375+2\cdot 2875=5375,
\ean
which also agrees with (\ref{gdata}).

From these results, we expect that ghost gravitational Gromov-Witten invariants can be derived from a 
combination of standard Gromov-Witten invariants and the ghost Gromov-Witten invariants $\langle{\cal O}_{h^{-1}}{\cal O}_{h^3}\rangle_{0,d}$ by applying the standard relations of gravitational Gromov-Witten invariants. This is because our proof of Theorem \ref{mainA} relies solely on these relations.

Then {\bf can we compute the ghost Gromov-Witten invariants $\langle{\cal O}_{h^{-1}}{\cal O}_{h^3}\rangle_{0,d}$ in a context other than mirror symmetry?} The answer is {\bf ``yes''} in a certain sense. Although we omit the proofs, we have the following residue integral representation of $\langle{\cal O}_{h^{a}}{\cal O}_{h^b}\rangle_{0,d}$ up to $d=3$, which is applicable to $M_{k}^{k}$ \cite{Jin2}:
\ban
\langle{\cal O}_{h^{a}}{\cal O}_{h^b}\rangle_{0,1}&=&-\frac{1}{2}\frac{1}{(2 \pi \sqrt{-1})^{2}}\oint_{C_{0}} \frac{dz_{0}}{(z_{0})^{k}} \oint_{C_{1}} \frac{dz_{1}}{(z_{1})^{k}}(z_{0}-z_{1})^2 e(z_{0},z_{1})
W_{a}(z_{0},z_{1})W_{b}(z_{0},z_{1}),
\ean
\ban
\langle{\cal O}_{h^{a}}{\cal O}_{h^b}\rangle_{0,2}&=&\frac{1}{2}\frac{1}{(2 \pi \sqrt{-1})^{3}}\oint_{C_{0}} \frac{dz_{0}}{(z_{0})^{k}} \oint_{C_{1}} \frac{dz_{1}}{(z_{1})^{k}}\oint_{C_{2}} \frac{dz_{2}}{(z_{2})^{k}}(z_{0}-z_{1})(z_{2}-z_{1}) \no\\
&&\times\frac{e(z_{0},z_{1})e(z_{1},z_{2})}{kz_{1}(2z_{1}-z_{0}-z_{2})} \left(W_{a}(z_{0},z_{1})+ W_{a}(z_{1},z_{2})\right) \left(W_{b}(z_{0},z_{1})+ W_{b}(z_{1},z_{2})\right),
\ean
\ban
\langle{\cal O}_{h^{a}}{\cal O}_{h^b}\rangle_{0,3}&=&\frac{1}{2}\frac{1}{(2 \pi \sqrt{-1})^{4}}\oint_{C_{0}} \frac{dz_{0}}{(z_{0})^{k}} \oint_{C_{1}} \frac{dz_{1}}{(z_{1})^{k}}\oint_{C_{2}} \frac{dz_{2}}{(z_{2})^{k}}\oint_{C_{3}} \frac{dz_{3}}{(z_{3})^{k}}(z_{0}-z_{1})(z_{3}-z_{2}) \no\\
&&\times\frac{e(z_{0},z_{1})e(z_{1},z_{2})e(z_{2},z_{3})}{kz_{1}(2z_{1}-z_{0}-z_{2})kz_{2}(2z_{2}-z_{1}-z_{3})} \left(W_{a}(z_{0},z_{1})+ W_{a}(z_{1},z_{2})+ W_{a}(z_{2},z_{3})\right)\no\\
&&\times \left(W_{b}(z_{0},z_{1})+ W_{b}(z_{1},z_{2})+ W_{b}(z_{2},z_{3})\right)
\ean
\ba
&&-\frac{1}{6}\frac{1}{(2 \pi \sqrt{-1})^{4}}\oint_{C_{(0)}} \frac{dz_{0}}{(z_{0})^{k}} \oint_{C_{(0)}} \frac{dz_{1}}{(z_{1})^{k}}\oint_{C_{(0)}} \frac{dz_{2}}{(z_{2})^{k}}\oint_{C_{(0)}} \frac{dz_{3}}{(z_{3})^{k}}\no\\
&&\times\frac{e(z_{0},z_{1})e(z_{0},z_{2})e(z_{0},z_{3})}{\left(kz_{0}\right)^2} \left(W_{a}(z_{0},z_{1})+ W_{a}(z_{0},z_{2})+ W_{a}(z_{0},z_{3})\right)\no\\
&&\times \left(W_{b}(z_{0},z_{1})+ W_{b}(z_{0},z_{2})+ W_{b}(z_{0},z_{3})\right),
\label{ggw}
\ea
where $\frac{1}{2}\frac{1}{2 \pi \sqrt{-1}}\oint_{C_{l}}dz_l$ has the same meaning as the one in Definition \ref{wvsdefi}, $\frac{1}{2}\frac{1}{2 \pi \sqrt{-1}}\oint_{C_{(0)}}dz_i$ means taking the residue at $z_{i}=0$,
and $W_{a}(z,w)$ is given by
\ban
W_{a}(z,w)=\frac{z^{a}-w^{a}}{z-w}.
\ean 
If we formally apply (\ref{ggw}) with $k=5$ and $\displaystyle{W_{-1}(z,w)=\frac{z^{-1}-w^{-1}}{z-w}=-\frac{1}{zw}}$ to $\langle{\cal O}_{h^{-1}}{\cal O}_{h^3}\rangle_{0,d}\;\;(d=1,2,3)$, we obtain
\ban
\langle{\cal O}_{h^{-1}}{\cal O}_{h^3}\rangle_{0,1}=-3250,\;\;
\langle{\cal O}_{h^{-1}}{\cal O}_{h^3}\rangle_{0,2}= -803125,\;\;
\langle{\cal O}_{h^{-1}}{\cal O}_{h^3}\rangle_{0,3}= -\frac{1686081250}{3}.\;\;
\ean
\end{ex}

\begin{ex} {\bf Ghost Gromov-Witten Invariants $\langle{\cal O}_{h^{k-2}}{\cal O}_{h^{-1}}\rangle_{0,d}$
of $M_{k}^{k}$}

In view of Theorem \ref{vmth} and Corollary \ref{fin}, the ghost Gromov-Witten invariant $\langle{\cal O}_{h^{k-2}}{\cal O}_{h^{-1}}\rangle_{0,d}$ is obtained by applying Theorem \ref{vmth} ``formally'' to the $j=0$ case:
\ban
\langle{\cal O}_{h^{k-2}}{\cal O}_{h^{-1}}\rangle_{0}(t)&=&w({\cal O}_{h^{k-2}}{\cal O}_{h^{-1}})_{0}(x(t)),
\ean  
where $w({\cal O}_{h^{k-2}}{\cal O}_{h^{-1}})_{0}(x)$ is given by
\ban
w({\cal O}_{h^{k-2}}{\cal O}_{h^{-1}})_{0}(x)&=&k\int^{x}\tilde{L}_{0}^{k,k}(u)du\\
                                                                     &=&k\int^{x}\left(1+\sum_{d=1}^{\infty}\frac{(kd)!}{(d!)^{k}}e^{du}\right)\\
                                                                     &=&kx+\sum_{d=1}^{\infty}\frac{k}{d}\cdot\frac{(kd)!}{(d!)^{k}}e^{dx}.
\ean
Therefore, we define the ghost Gromov-Witten invariant $\langle{\cal O}_{h^{k-2}}{\cal O}_{h^{-1}}\rangle_{0,d}$ as follows.
\begin{defi}
Let $M_{k}^{k}$ be a smooth degree $k$ hypersurface in $CP^{k-1}$ and let $t=t(x)=\frac{W_{(1)}^{k,k}(x)}{\tilde{W}_{(0)}^{k,k}(x)}$ be the mirror map used in the mirror computation of Gromov-Witten invariants of $M_{k}^{k}$. We define the ghost Gromov-Witten invariant $\langle{\cal O}_{h^{k-2}}{\cal O}_{h^{-1}}\rangle_{0,d}$ as follows:
\ban
kt+\sum_{d=1}^{\infty}\langle{\cal O}_{h^{k-2}}{\cal O}_{h^{-1}}\rangle_{0,d}e^{dt}
:=kx(t)+\sum_{d=1}^{\infty}\frac{k}{d}\cdot\frac{(kd)!}{(d!)^{k}}e^{dx(t)},
\ean
where $x=x(t)$ is the inverse of the mirror map.
\end{defi}
\end{ex}

From this definition, we can compute $\langle{\cal O}_{h^{k-2}}{\cal O}_{h^{-1}}\rangle_{0,d}$
to arbitrary order. 
In Appendix A, we present tables of these invariants for $k=3$ to $k=5$ up to $d=10$. Note that $\langle{\cal O}_{h^{k-2}}{\cal O}_{h^{-1}}\rangle_{0,d}$ is non-vanishing in the $k=3,4$ cases, i.e., the cases of the elliptic curve and the K3 surface, where no rational curves of any degree exist. Looking at the numerical data,
we assume the multiple cover formula for these invariants and consider the following expansion: 
\ban
kt+\sum_{d=1}^{\infty}\langle{\cal O}_{h^{k-2}}{\cal O}_{h^{-1}}\rangle_{0,d}e^{dt}
=kt+\sum_{d=1}^{\infty}n^{k,k}_{gh.}(d)\left(-\log(1-e^{dt})\right).
\ean
We also present $n^{k,k}_{gh.}(d)$ in the tables. 

Lastly, we present a conjecture on $n^{k,k}_{gh.}(d)$
obtained from our numerical data.
\begin{conj}
$n^{k,k}_{gh.}(d)$ is always an integer and $\frac{n^{3,3}_{gh.}(d)}{d}$ is also always an integer. Moreover,
$\frac{n^{3,3}_{gh.}(d)}{d}$ is given by 
\ban
\frac{n^{3,3}_{gh.}(d)}{d}=\left\{\begin{array}{rr}-27&(d=3i-2)\\27&(d=3i-1)\\0&(d=3i)\end{array}
\right. .
\ean
\end{conj}

\newpage

\appendix

\section{Tables of Ghost Gromov-Witten Invariants $\langle{\cal O}_{h^{k-2}}{\cal O}_{h^{-1}}\rangle_{0,d}$}

\begin{table}[H]
\centering
\caption{$k=3$} \label{k=3}

\begin{tabular}{|c|r|r|r|r|}
\hline
$d$&$\langle{\cal O}_{h}{\cal O}_{h^{-1}}\rangle_{0,d}$& $n^{3,3}_{gh}(d)$ & $\frac{n^{3,3}_{gh}(d)}{d}$  \\
\hline
1&$-27$ &$-27$&$-27$ \\ \hline
2&$81/2$ &$54$& $27$\\ \hline
3&$-9$&$0$&$0$ \\ \hline
4&$-351/4$&$-108$&$-27$ \\ \hline
5&$648/5$ & $135$&$27$\\ \hline
6&$27/2$ & $0$&$0$\\ \hline
7&$-1350/7$ &$-189$&$-27$\\ \hline
8&$1377/8$ & $216$& $27$\\ \hline
9&$-3$& $0$&$0$\\ \hline
10&$-972/5$& $-270$&$-27$\\ \hline
\end{tabular}
\end{table}

 \begin{table}[H]
\centering
\caption{$k=4$} \label{k=4}

\begin{tabular}{|c|r|r|r|r|}
\hline
$d$&$\langle{\cal O}_{h^2}{\cal O}_{h^{-1}}\rangle_{0,d}$& $n^{4,4}_{gh}(d)$   \\
\hline
1&$-320$ &$-320$ \\ \hline
2&$-800$ &$-640$\\ \hline
3&$-8960/3$&$-2880$\\ \hline
4&$-2960$&$-2560$\\ \hline
5&$-8064$ & $-8000$\\ \hline
6&$-22400/3$ & $-5760$\\ \hline
7&$-110080/7$ &$-15680$\\ \hline
8&$-11720$ & $-10240$\\ \hline
9&$-242240/9$& $-25920$\\ \hline
10&$-20160$&$-16000$\\ \hline
\end{tabular}
\end{table}
 \begin{table}[H]
\centering
\caption{$k=5$} \label{k=5}

\begin{tabular}{|c|r|r|r|r|}
\hline
$d$&$\langle{\cal O}_{h^3}{\cal O}_{h^{-1}}\rangle_{0,d}$& $n^{5,5}_{gh}(d)$   \\
\hline
1&$-3250$ &$-3250$ \\ \hline
2&$-803125$ &$-801500$\\ \hline
3&$-1686081250/3$&$-562026000$\\ \hline
4&$-1089994203125/2$&$-544996700000$\\ \hline
5&$-626257527494400$ & $-626257527493750$\\ \hline
6&$-2396497053884609125/3$ & $-798832351013589000$\\ \hline
7&$-7656135988083727275000/7$ &$-1093733712583389610250$\\ \hline
8&$-6304748148021421341403125/4$ & $-1576187037005082836800000$\\ \hline
9&$-21251095656353576496361281250/9$& $-2361232850705952943852800000$\\ \hline
10&$-3646338072561587110650897395000$&$-3646338072561586797522133487500$\\ \hline
\end{tabular}
\end{table}

\end{document}